\documentclass[pdftex,12pt,a4paper]{article} 

\usepackage{verbatim}
\usepackage{amsbsy}
\usepackage{amscd}
\usepackage[mathscr]{eucal}
\usepackage{amsmath,amsthm, amssymb, url, enumerate}
\usepackage[pdftex]{graphicx} 
\usepackage{cases}
\usepackage{stmaryrd}
\usepackage[T1]{fontenc}
\usepackage{textcomp}
\usepackage{color}
\usepackage[all]{xy}
\usepackage{cases}
\usepackage{stmaryrd}
\usepackage{bm}
\usepackage{here}
\usepackage{comment}

\usepackage[svgnames]{xcolor}%
\usepackage{tikz}
\usepackage{tikz-cd}
 \usetikzlibrary{calc}
\usetikzlibrary {arrows.meta}
\usetikzlibrary{arrows}

\definecolor{foBlue}{RGB}{35,105,153}
\definecolor{foOrange}{RGB}{155,77,18}
\definecolor{foRed}{RGB}{181,36,52}
\definecolor{foMesh}{RGB}{171,179,189}

\newtheorem{defi}{Definition}[section]
\newtheorem{thm}[defi]{Theorem}
\newtheorem{prop}[defi]{Proposition} 
\newtheorem{lem}[defi]{Lemma}
\newtheorem{rem}[defi]{Remark}
\newtheorem{cor}[defi]{Corollary}

\newtheorem{fact}[defi]{Fact}
\newtheorem{ex}[defi]{Example}
\newenvironment{prf}
  {\par\noindent\textit{Proof.}\ \normalfont\ignorespaces}
  {\unskip\nobreak\hfill$\square$\par\medskip}

\newcommand{\compl}[1]{\mathbb{C}^{#1}}
\newcommand{\real}[1]{\mathbb{R}^{#1}}
\newcommand{\ratio}[1]{\mathbb{Q}^{#1}}
\newcommand{\inte}[1]{\mathbb{Z}^{#1}}
\newcommand{\natu}[1]{\mathbb{N}^{#1}}

\newcommand{\age}{\mathrm{age}}

\newcommand{\Hilb}{\mathrm{Hilb}}

\begin{document}

\title{Euler characteristics and moderateness of
Fujiki--Oka resolutions}
\author{Yusuke\ Sato}
\date{}
\maketitle 

%











\footnote{\hspace{-0.7cm} \textbf{2020 Mathematics Subject Classification:}
Primary 14M25  : Secondary 14E16. 14J17 }


\maketitle

\begin{abstract}
Fujiki--Oka resolutions of cyclic quotient singularities
of type \(\frac1r(1,a_2,\ldots,a_n)\) are constructed using
multidimensional continued fractions.
In this paper, we compute the Euler characteristics of
Fujiki--Oka resolutions in arbitrary dimension using the
heights of the coefficients of the associated remainder
polynomials. We also give a necessary and sufficient condition for moderateness of Fujiki--Oka resolutions in dimension three.


\end{abstract}

\section{Introduction}

\subsection{Background}

The classical McKay correspondence gives a bijection between the nontrivial irreducible representations of a finite subgroup \(G\subset\mathrm{SL}(2,\mathbb C)\) and the irreducible exceptional curves in the minimal resolution of \(\mathbb C^2/G\).
In higher dimensions, crepant resolutions of
\(\mathbb C^n/G\), with \(G\subset\mathrm{SL}(n,\mathbb C)\),
provide a geometric setting for extensions of this correspondence,
when such resolutions exist.
The following theorem of Batyrev gives a numerical form of
the McKay correspondence.

\begin{thm}\upshape{\cite{Batyrev}} \label{Batyrev}
If a projective crepant resolution $Y \to \compl{n}/G$ exists, then the Euler characteristic of $Y$ equals the number of conjugacy classes of $G$. 
\end{thm}

For cyclic groups, this number is $|G|$. To study analogous relations for cyclic quotient singularities beyond the Gorenstein case, it is useful to consider other distinguished classes of
toric resolutions.

Let $\sigma$ be a full-dimensional strongly convex rational
polyhedral cone in a lattice $N$. A Hilbert basis resolution of
$X_\sigma$ is a toric resolution whose primitive ray generators
are exactly the elements of the Hilbert basis of $\sigma \cap N$.
Every three-dimensional normal affine toric variety admits
a Hilbert basis resolution
\cite[Theorem~2.9]{BGS}.
A smooth subdivision of $\sigma$ is called moderate if the affine hyperplane through its primitive ray generators meets every extremal ray of $\sigma$ for each
maximal cone.
For toric resolutions of a $\mathbb{Q}$-Gorenstein affine toric
variety, one has
\[
  \text{crepant}
  \;\Longrightarrow\;
  \text{moderate}
  \;\Longrightarrow\;
  \text{Hilbert basis}.
\]
 Moderate toric resolutions were introduced by
Ch\'avez-Mart\'inez, Duarte, and Yasuda in their study of
$F$-blowups and essential divisors \cite{CMDY}.
The relevant essential divisors, in the sense of Bouvier and
Gonzalez-Sprinberg \cite{BGS}, are the torus-invariant divisors over
$X_\sigma$ whose centers are prime divisors on every toric
resolution; they correspond to the elements of the Hilbert
basis of $\sigma\cap N$.
The existence of a moderate toric resolution ensures that
the normalized limit $F$-blowup contains prime divisors
representing all of these essential divisors.
In dimension three, Kaijima and Yamamoto show that every
moderate toric resolution is a minimal terminal
\(\mathbb Q\)-factorial model.
In particular, it is crepant when the original affine toric
variety is canonical \cite{KY}.

Fujiki--Oka resolutions provide an explicit class of toric
resolutions of cyclic quotient singularities of type
$\frac{1}{r}(1,a_2,\ldots,a_n)$.
Their construction is described by multidimensional continued
fractions and the associated remainder polynomials
\cite{Ashikaga}.
In dimension two, they coincide with the minimal resolutions,
while in the three-dimensional Gorenstein case they are crepant
\cite{SS}. Hilbert basis properties of Fujiki--Oka resolutions and iterated
Fujiki--Oka resolutions for three-dimensional canonical cyclic
quotient singularities were studied in \cite{SS2}.
Our aim is to characterize the moderateness of the Fujiki--Oka
resolution from its continued-fraction data.
In dimension three, we prove that this resolution is moderate
if and only if every coefficient of its remainder
polynomial has nonpositive height.
We also express its Euler characteristic in arbitrary dimension
as the group order plus the sum of the coefficient heights.

\subsection{Main results}


Let $G$ be a finite cyclic subgroup of ${\rm GL}(n,\compl{})$ generated by ${\rm diag}(\varepsilon^{a_1},\dots,\varepsilon^{a_n})$ where $\varepsilon$ is a primitive $r$-th root of unity. For simplicity, ${\rm diag}(\varepsilon^{a_1},\dots,\varepsilon^{a_n})$ is denoted by $\frac{1}{r}(a_1,\dots,a_n)$.\\
Let \(N=\mathbb Z^n\) with standard basis \(\bm{e}_1,\ldots,\bm{e}_n\), and let $N_G$ be $\inte{n}+\inte{}\cdot \frac{1}{r}(a_1,\dots,a_n)$.
We set a rational strongly convex polyhedral cone $\sigma \subset N_{G}\otimes_{\inte{}}\real{}$ as $\real{}_{\geq0}\bm{v}_1 + \cdots + \real{}_{\geq0}\bm{v}_n$ where $\real{}_{\geq0}$ is the set of all non negative elements in $\real{}$ and $\bm{v_i}$ are in $N$. For simplicity, $\real{}_{\geq0}\bm{v}_1 + \cdots + \real{}_{\geq0}\bm{v}_n$ is denoted by ${\rm Cone}(\bm{v}_1,\dots,\bm{v}_n)$.
The {\it dimension} of a cone $\sigma$ is defined as the dimension of $\real{}\cdot\sigma$ as a vector space over $\real{}$. 
To shorten notation, $\sigma$ also signifies the finite fan consisting of all faces of $\sigma$. 
Let $\sigma_0$ be a cone generated by $\bm{e}_1,\ldots, \bm{e}_n$, that is ${\rm Cone}(\bm{e}_1,\ldots, \bm{e}_n)$. Then the toric variety $X(N_G,\sigma_0)$ determined by $N_G$ and the finite fan $\sigma_0$ is $\mathbb{C}^n/G$.

\begin{defi}\upshape
Let $\sigma=\operatorname{Cone}(\bm{v}_1,\ldots,\bm{v}_n)$ be an
$n$-dimensional cone, and let $\bm{v}=(a_1,\ldots,a_n)/r$ in $(\sigma\cap N_G)\setminus\{0\}$. Then \(\Sigma_{\bm{v}}\) is the fan consisting of the
\(n\)-dimensional cones among \(\sigma_1,\ldots,\sigma_n\)
and all their faces, where $ \sigma_i
 =\operatorname{Cone}
  (\bm{v}_1,\ldots,\bm{v}_{i-1},\bm{v}_{i+1},\ldots,\bm{v}_n,\bm{v})$ with
 $i=1,\ldots,n$.
The subdivision $\Sigma_{\bm{v}} \to \sigma$ is called the star subdivision. 
In addition, we call the induced toric morphism $f:X(N_G,\Sigma_{\bm{v}}) \to X(N_G,\sigma)$ "\textbf{the weighted blow-up}" of $X(N_G,\sigma)$ with weight $\bm{v}=\frac{1}{r}(a_1,\dots,a_n)$.
\end{defi}

\begin{defi}\upshape
Let \(c=(b_1,\ldots,b_n)/r\) be an \(n\)-dimensional fraction,
where \(r\) is a positive integer and \(b_1,\ldots,b_n\) are
integers satisfying \(0\le b_i<r\) for \(1\le i\le n\).
We define its height by
\[
 h(c):=\sum_{i=1}^n b_i-r.
\]
\end{defi}

We write \(h(\mathcal R_*(v))\) for the sum of the heights
of its nonzero coefficient occurrences (see Section~3).

\begin{thm}\label{intro:Eular}
Let \(G\) be a cyclic group of type
\(\frac{1}{r}(1,a_2,\ldots,a_n)\), with $\gcd(r,a_2,\dots,a_n)=1$.
For the Fujiki--Oka resolution
\(f:Y\to\compl{n}/G\), we have
\[
 \chi(Y)
 =
 h\left(
   \mathcal R_*\left(\frac{(1,a_2,\ldots,a_n)}{r}\right)
 \right)+r,
\]
where \(\chi(Y)\) is the topological Euler characteristic of \(Y\).
\end{thm}

The definition of a Fujiki--Oka resolution and a remainder polynomial $\mathcal{R}_*$ are introduced in section 2. 
Since \(|G|=r\), Theorem~\ref{intro:Eular} expresses the difference \(\chi(Y)-|G|\) as the sum of the coefficient heights.
Therefore, it provides a numerical analogue of the McKay correspondence for this class of cyclic quotient singularities.
Our second main result gives a criterion for moderateness
in dimension three.

\begin{thm}
\label{intro:moderate}
\upshape
Let \(G\) be of type \(\frac1r(1,a_2,a_3)\) with \(\gcd(r,a_2,a_3)=1\).
The associated Fujiki--Oka resolution is moderate if and
only if every nonzero coefficient of $ \mathcal R_*\left(\frac{(1,a_2,a_3)}r\right)$
has nonpositive height.
\end{thm}

Under these equivalent conditions, the resolution is a
Hilbert basis resolution and satisfies $ \chi(Y)\le |G|$,
with equality if and only if it is crepant. The Hilbert basis property follows directly from moderateness.
The Euler characteristic formula involves the sum of the
coefficient heights, whereas the moderateness criterion
requires each coefficient height to be nonpositive.
Without the moderateness assumption, the equality \(\chi(Y)=|G|\) does not imply
the Fujiki-Oka resolution is crepant. (see Example~3.9).

\par
This paper is organized as follows.
Section 2 introduces the resolution of singularities (that is called the Fujiki--Oka resolution) determined by the multidimensional continued fraction.
In section 3, we show that the Euler characteristic of the Fujiki--Oka resolution is calculated from the continued fraction.
In section 4, we establish the criterion for moderateness of Fujiki--Oka resolutions in dimension three. Section 5 presents several examples related to the main results.\\

\textbf{Acknowledgments}:
The author is grateful to Professor Yujiro Kawamata, whose question at the author's doctoral dissertation defense prompted this study of coefficient heights. The author also would like to extend gratitude to Professor Tadashi Ashikaga, Professor Yukari Ito, and Associate Professor Kohei Sato for their several helpful comments.


\section{Toric resolution via multidimensional continued fraction}
In this section, we introduce the Fujiki--Oka resolution as a toric resolution determined by a multidimensional continued fraction. The Fujiki--Oka construction originates in the work of Fujiki \cite{Fujiki} and Oka \cite{Oka}. We follow Ashikaga \cite{Ashikaga} for its description
in terms of multidimensional continued fractions.
The multidimensional continued fraction consists of the remainder polynomial and the round down polynomial.
As described later, the remainder polynomial gives the type of quotient singularities appearing at each stage
of a Fujiki--Oka resolution. On the other hand, the round down polynomial gives a part of the intersection numbers of exceptional curves.
As the round down polynomial is not necessary for the discussion in this paper, we introduce only the definition of the remainder polynomial.

\subsection{Remainder polynomial}\label{ACF}
\begin{defi}\upshape\label{profrac}
Let $n$ and $r$ be elements of $\natu{}$. Let $\mathbf{a}=(a_1,\dots,a_n) \in \inte n$ which satisfies $0\leq  a_i \leq r-1$ for $1\leq i \leq n$. We call the symbol
$$
\frac{\mathbf{a}}{r}=\frac{(a_1,\dots,a_n)}{r}
$$
an {\it $n$-dimensional proper fraction}.\\
In addition, a proper fraction is called {\it semi-unimodular} if at least one of its numerator entries is $1$.

\end{defi}

In the following, the symbol $\ratio{prop}_n$ (resp. $\overline{\ratio{prop}_n}$) means the set of $n$-dimensional proper fractions (resp. the set $\ratio{prop}_n \cup \{\infty\}$). Moreover, $\overline{\ratio{prop}_n}[x_1,\dots,x_n]$  denotes the set consisting of all noncommutative polynomials with $n$ variables over $\overline{\ratio{prop}_n}$. 
The remainder polynomial is obtained via  {\it remainder maps} for a {\it semi-unimodular proper fraction}.

\begin{defi}\upshape (\cite[Def 3.1.]{Ashikaga})
Let $\frac{(a_1,a_2,\dots,a_n)}{r}$ be a semi-unimodular proper fraction.  For $1\leq i \leq n$, the {\it $i$-th remainder map} $R_i:\overline{\ratio{prop}_n} \to \overline{\ratio{prop}_n}$ is defined by
$$
R_i\left(\frac{(a_1, a_2,\dots,a_n)}{r}\right)=
\left\{ 
\begin{array}{cc}
 \frac{\left(\overline{a_1}^{a_i},\ \dots,\ \overline{-r}^{a_i},\ \dots, \overline{a_n}^{a_i}\right)}{a_i}
         & {\rm if} \  a_i \neq 0\\  
 \infty &  {\rm if}\  a_i=0
\end{array}
\right.
$$
and $R_i(\infty)=\infty$ where $\overline{a_j}^{a_i}$ is an integer satisfying $0\leq \overline{a_j}^{a_i} < a_i$ and $\overline{a_j}^{a_i} \equiv a_j$ modulo $a_i$.
\end{defi}

\begin{ex}\upshape
If $v=\frac{(1,2,7)}{12}$, then we have
$$
R_2(v)=\frac{(1,0,1)}{2}\ \text{and} \ R_3(v)=\frac{(1,2,2)}{7}.
$$

\end{ex}

\begin{defi}\upshape\cite[Def 3.2.]{Ashikaga}\label{DOACF}
Let $\frac{\mathbf{a}}{r}$ be an $n$-dimensional semi-unimodular proper fraction and $I=\{1,\dots,n \}$. The {\it remainder polynomial} $\mathcal{R}_*\left(\frac{\mathbf{a}}{r}\right) \in \overline{\ratio{prop}_n}[x_1,x_2,\dots,x_n]$ is defined by
  $$
  \mathcal{R}_*\left(\frac{\mathbf{a}}{r}\right)=\frac{\mathbf{a}}{r}+
                                      \sum_{(i_1,i_2,\dots,i_l)\in \mathbf{I}^l,\: l\geq 1 }(R_{i_l}\cdots R_{i_2}R_{i_1})\left(\frac{\mathbf{a}}{r}\right)\cdot x_{i_1}x_{i_2}\cdots x_{i_l},
  $$
  where we exclude terms with coefficient $\infty$ or $\frac{(0,0,\dots,0)}{1}$.
 
\end{defi}

The variables $x_1,\ldots,x_n$ are noncommuting.
The monomial $x_{i_1}\cdots x_{i_l}$ records the order
in which the remainder maps are applied; $R_{i_1}$ is applied first, followed successively by
$R_{i_2},\ldots,R_{i_l}$. We note that distinct monomials are kept as separate terms even when their coefficients coincide. 
Since the denominator of the image of a remainder map is monotonically decreasing, the remainder polynomial consists of a finite number of terms.

\begin{ex}\label{remainderex}\upshape
Let $v=\frac{(1,2,7)}{12}$, then the remainder polynomial is
\begin{eqnarray*}
\mathcal{R}_*\left(\frac{(1,2,7)}{12}\right)&=& \frac{(1,2,7)}{12}+\frac{(1,0,1)}{2}x_2+\frac{(1,2,2)}{7}x_3\\ 
                   &+&\frac{(1,1,0)}{2}x_3x_2+\frac{(1,0,1)}{2}x_3x_3.\\ 
\end{eqnarray*}

\end{ex}

\subsection{Fujiki--Oka resolution}
In this subsection, we introduce the Fujiki--Oka resolution as a toric resolution. We fix a primitive lattice point $\bm{v}=\frac{1}{r}(a_1,\dots,a_n)$ in $N_G$. Assume that $\bm{v}$ generates $N_G/N$. We denote by $\sigma_k$ the $n$-dimensional cone ${\rm Cone}(\bm{e}_1,\dots,\bm{e}_{k-1},\bm{v}, \bm{e}_{k+1},\dots,\bm{e}_n)$ for \(k\) with \(a_k>0\).
We consider the star subdivision $\Sigma_{\bm{v}} \to \sigma={\rm Cone}(\bm{e}_1,\dots,\bm{e}_n)$ at $\bm{v}$, where the fan $\Sigma_{\bm{v}}$ consists of the cones $\sigma_k$ and all their faces. 
Let $N_{G,k}$ be the sublattice of $N_G$ which is generated by $\bm{e}_1, \dots,  \bm{e}_{k-1}, \bm{v},$ $\bm{e}_{k+1}, \dots , \bm{e}_n$. Then the dual lattice $M_{G,k}:={\rm Hom}(N_{G,k},\inte{})$ has dual basis $\{\xi_1,\dots, \xi_n\}$ which satisfy:
$$
\xi_j=\begin{cases}
 x_jx_k^{-\frac{a_j}{a_k}}
         & {\rm if} \  j \neq k,\\  
 x_k^{\frac{r}{a_k}} &  {\rm if}\  j=k.
\end{cases}
$$
\noindent
Therefore, the affine toric variety \(X(N_G,\sigma_k)\)
is a cyclic quotient of type\\ \(\frac{1}{a_k}(a_1,\ldots,a_{k-1},-r,a_{k+1},\ldots,a_n)\), where \(\xi_1,\ldots,\xi_n\) are coordinates on the covering affine space \(X(N_{G,k},\sigma_k)\cong\mathbb C^n\). Let $G=\frac{1}{r}(1,a_2,\dots,a_n)$ and $\bm{v}=\frac{1}{r}(1,a_2,\dots,a_n)$ in $N_G$ with $\gcd(r,a_2,\dots, a_n)=1$. Then $\bm{v}$ generates $N_G/N$. 
The primitive lattice point \(\bm{v}\) is called the \emph{Oka center} of \(\sigma\) with respect to the distinguished generator \(\bm{e}_1\).
More generally, let
\(\tau=\operatorname{Cone}(\bm{e}_1,\bm{w}_2,\ldots,\bm{w}_n)\)
be a singular cone occurring in the construction, and let
\((1,b_2,\ldots,b_n)/s\)
be its local proper fraction with respect to these ordered generators.
Its Oka center is the primitive lattice point
\[
\bm{u}=\frac{{\bm e}_1+b_2\bm{w}_2+\cdots+b_n\bm{w}_n}{s}\in N_G.
\]
Thus, the proper fraction records the coordinates of the Oka center in the ordered real basis of the local cone. Throughout the construction, all Oka centers are regarded as lattice points in the fixed lattice \(N_G\).

We consider the star subdivision at $\bm{v}$. Then the types of quotient singularities appearing in this subdivision equal  the image of the $k$-th remainder map $R_k\left((1,a_2,\dots,a_n)/r\right)$.
Whenever $R_k((1,a_2,\ldots,a_n)/r)$ is a nonzero proper fraction, it gives the coordinates of the Oka center of $\sigma_k$ with respect to the ordered generators $(\bm{e}_1,\ldots,\bm{e}_{k-1},\bm{v},\bm{e}_{k+1},\ldots,\bm{e}_n)$.
We repeat the star subdivision at these Oka centers. The types of quotient singularities  appearing at each stage of star subdivisions are obtained from the remainder polynomial $\mathcal{R}_*\left((1,a_2,\dots,a_n)/r\right)$. 
Thus, the remainder polynomial induces the toric morphism which is called the Fujiki--Oka resolution (see \cite{Ashikaga} for more details).

\begin{ex}\upshape
Let \(G=\frac{1}{12}(1,2,7)\), as in Example~\ref{remainderex}.
The star subdivision at \(\bm{v}_1=\frac{1}{12}(1,2,7)\) has three maximal
cones, \(\sigma_1,\sigma_2,\sigma_3\), shown on the right of Figure~$1$.
The coefficients of the remainder polynomial indicate that the affine toric variety corresponding to $\sigma_2$ (resp. $\sigma_3$) has a quotient singularity of type $\frac{1}{2}(1,0,1)$ (resp. $\frac{1}{7}(1,2,2)$). Since the image of the remainder map $R_1((1,2,7)/12)$ is $(0,0,0)/1$, the affine toric variety $X(N_G,\sigma_1)$ is smooth.\par
Next, we consider the star subdivision of \(\sigma_3\)
at the point with coordinates \(\frac{1}{7}(1,2,2)\)
relative to \((\bm{e}_1,\bm{e}_2,\bm{v}_1)\), that is, at
\(\bm{v}_2=\frac{1}{12}(2,4,2)\in N_G\).
 After that, the type of quotient singularity corresponding to $\sigma_{32}$ is $\frac{1}{2}(1,1,0)$, and $\sigma_{33}$ is $\frac{1}{2}(1,0,1)$. Finally, we subdivide at $\bm{v}_6=\frac{1}{12}(6,0,6)$ and   $\bm{v}_7=\frac{1}{12}(7,2,1)$, then we have the subdivision $\Sigma \to \sigma$ where $\Sigma$ consists of eight three-dimensional cones and their faces (the left of Figure 1). The monomials \(x_2\) and \(x_3^2\) have the same coefficient
\((1,0,1)/2\), but they give different Oka centers,
\(\bm{v}_6\) and \(\bm{v}_7\), respectively. Conversely, \(x_3x_2\) and \(x_3^2\) have different coefficients but give the same Oka center
\(\bm{v}_7=(\bm{e}_1+\bm{v}_2)/2\).
This subdivision induces the toric morphism $f:X(N_G,\Sigma) \to \compl{3}/G$, and it is the Fujiki--Oka resolution.
From now on, $\Sigma_{ {\rm F.O.R.}}$ denotes the fan which gives the Fujiki--Oka resolution.

\begin{figure}[htbp]
\centering
    \begin{tabular}{ccc}
      \begin{minipage}{0.3\columnwidth}
  
\begin{tikzpicture}
\draw (0,0) -- (0:4) -- ++(120:4)  -- cycle;
\coordinate (6) at (60:2) node at (6) [left] {$\bm{v}_6$};
\coordinate (A) at (0,0) ;
\coordinate (B) at (0:4) ;
\coordinate (C) at (60:4) ;
\coordinate (2) at ($(B)!.5!(6)$) node at (2) [below] {$\bm{v}_2$};
\coordinate (7) at ($(C)!.5!(2)$)  node at (7) [left] {$\bm{v}_7$};
\coordinate (1) at ($(A)!.5!(2)$)  ;
\draw (A) -- (1);
\draw (B) -- (1);
\draw (C) -- (1);
\draw (1) -- (6);
\draw (1) -- (2);
\draw (1) -- (7);
\draw (B) -- (2);
\draw (B) -- (7);
\draw (C) -- (2);
\coordinate (P) at (4,2) ;  
\coordinate (Q) at (5.5,2) ;

\draw [->] (P) -- (Q);

\fill (C) circle (2pt)  (A) circle (2pt) (B) circle (2pt)  (6) circle (2pt) (2) circle (2pt) (1) circle (2pt) (7) circle (2pt); 

      \end{tikzpicture}

       \end{minipage} &
      \begin{minipage}{0.3\columnwidth}
  
\begin{tikzpicture}
\draw (0,0) -- (0:4) -- ++(120:4)  -- cycle;
\coordinate (6) at (60:2) ;
\coordinate (A) at (0,0) ;
\coordinate (B) at (0:4) ;
\coordinate (C) at (60:4) ;
\coordinate (2) at ($(B)!.5!(6)$)  node at (2) [below] {$\bm{v}_2$} ;
\coordinate (7) at ($(C)!.5!(2)$)  ;
\coordinate (1) at ($(A)!.5!(2)$) node at (1) [above left=0.01mm] {$\bm{v}_1$};
\draw (A) -- (1);
\draw (B) -- (1);
\draw (C) -- (1);
\draw (1) -- (2);
\draw (B) -- (2);
\draw (C) -- (2);

\fill (C) circle (2pt)  (A) circle (2pt) (B) circle (2pt)  (6) circle (2pt) (2) circle (2pt) (1) circle (2pt) (7) circle (2pt); 

\coordinate (P) at (4,2) ;
\coordinate (Q) at (5.5,2) ;

\draw [->] (P) -- (Q);

       \coordinate (sigma32) at ($(2)!.4!(6)$) node at (sigma32)  {$\sigma_{32}$};
             \coordinate (sigma33) at ($(7)!.4!(B)$) node at (sigma33)  {$\sigma_{33}$};
      \end{tikzpicture}

      \end{minipage}
      \begin{minipage}{0.3\columnwidth}
  
 \begin{tikzpicture}
\draw (0,0) -- (0:4) -- ++(120:4)  -- cycle;
\coordinate (6) at (60:2) ;
\coordinate (A) at (0,0) node at (A) [left]{$\bm{e}_3$} ;
\coordinate (B) at (0:4) node at (B) [right]{$\bm{e}_2$};
\coordinate (C) at (60:4) node at (C) [above]{$\bm{e}_1$};
\coordinate (2) at ($(B)!.5!(6)$)  ;
\coordinate (7) at ($(C)!.5!(2)$)  ;
\coordinate (1) at ($(A)!.5!(2)$) node at (1) [above right=0.001mm] {$\bm{v}_1$};
\draw (A) -- (1);
\draw (B) -- (1);
\draw (C) -- (1);

\fill (C) circle (2pt)  (A) circle (2pt) (B) circle (2pt)  (6) circle (2pt) (2) circle (2pt) (1) circle (2pt) (7) circle (2pt); 

\coordinate (sigma1) at ($(A)!.5!(2)$) node at (sigma1) [below] {$\sigma_1$};
       \coordinate (sigma2) at ($(1)!.5!(6)$) node at (sigma2) [above=0.1mm] {$\sigma_2$};
       \coordinate (sigma3) at ($(2)!.5!(7)$) node at (sigma3)  {$\sigma_3$};
      \end{tikzpicture}

      \end{minipage}
     
    \end{tabular}
    \caption{Fujiki--Oka resolution for $G=\frac{1}{12}(1,2,7)$}
\end{figure}

\end{ex}

\section{Euler characteristic of Fujiki--Oka resolution}
In this section, we prove an Euler characteristic formula for Fujiki–Oka resolutions.

\begin{defi}\upshape
Let $\mathcal R_*(v)=\sum_X c_X X$,
where \(X\) runs through the monomials with nonzero coefficients.
We define the height of the remainder polynomial by
\[
 h(\mathcal R_*(v)):=\sum_X h(c_X).
\]
\end{defi}

\begin{defi}\upshape
The size of the remainder polynomial $S\left(\mathcal{R}_*\left(\frac{(a_1,\dots,a_n)}{r}\right)\right)$ is defined by the total number of numerator entries equal to \(1\), counted over all monomials.
\end{defi}

\begin{ex}\upshape

Let $c=(1,2,7)/12$, then ${\rm h}\left(c\right)=1+2+7-12=-2$.
By Example 2.5, ${\rm h}\left(\mathcal{R}_*(c)\right)=-2+0+-2+0+0=-4$, and $S(\mathcal{R}_*(c))=8$.

\end{ex}


The image of the remainder map in the component whose numerator is $1$ is $(0,\dots,0)/1$. It corresponds to the smooth cone of maximal dimension in $\Sigma_{ {\rm F.O.R.}}$. It follows that ${\rm S}(\mathcal{R}_*)$ equals the number of the cones of maximal dimension in $\Sigma_{{\rm F.O.R.}}$.

\begin{prop}\label{main1prop}
Let $c$ be an $n$-dimensional semi-unimodular proper fraction $(a_1,\dots,a_n)/r$, then the following holds:
$$
{\rm S}(\mathcal{R}_*(c))={\rm h}(\mathcal{R}_*(c)) + r.
$$

\end{prop}

\begin{proof}
Let \(k\) be the maximum of the numerators of \(c\).
The proof is by induction on \(k\).
We first show the case \(k=1\).
By permuting the components, we may assume that $a_1=\cdots=a_s=1, a_{s+1}=\cdots=a_n=0$ for some positive integer \(s\).
Since $ \mathcal R_*(c)=(1,\ldots,1,0,\ldots,0)/r$
we have $S(\mathcal R_*(c))
 =s=s-r+r=h(\mathcal R_*(c))+r$.

Next, assume that the proposition holds when all numerators
are at most \(k\). We prove it for \(k+1\).
By permuting the components, we may write $c=(a_1,\ldots,a_s,b_{s+1},\ldots,b_t,1,\ldots,1)/r$,
where \(a_i=k+1\) for \(1\le i\le s\), and
\(b_j=0\) or \(2\le b_j\le k\) for \(s+1\le j\le t\).

For \(i=1,\ldots,s\), the image of the \(i\)-th remainder map is
\[
 R_i(c)=
 \frac{
 (0,\ldots,0,\overline{-r}^{\,k+1},0,\ldots,0,
  \overline{b_{s+1}}^{\,k+1},\ldots,
  \overline{b_t}^{\,k+1},1,\ldots,1)
 }{k+1}.
\]
All its numerators are at most \(k\), so the induction
hypothesis applies to \(\mathcal R_*(R_i(c))\).
Similarly, it applies to \(\mathcal R_*(R_j(c))\)
whenever \(b_j>1\). In the sums below, expressions involving
\(\mathcal R_*(R_i(c))\) are omitted whenever
\(R_i(c)=\infty\) or \(R_i(c)=(0,\ldots,0)/1\).
By the above argument, we have the following formula:
\begin{align*}
 S(\mathcal R_*(c))
 &=\sum_{i=1}^n S(\mathcal R_*(R_i(c)))+n-t\\
 &=\sum_{i=1}^n h(\mathcal R_*(R_i(c)))
   +\sum_{i=1}^s a_i+\sum_{i=s+1}^t b_i+n-t\\
 &=\sum_{i=1}^n h(\mathcal R_*(R_i(c)))+h(c)+r\\
 &=h(\mathcal R_*(c))+r.
\end{align*}
Thus, by induction, the proposition holds for all
\(n\)-dimensional semi-unimodular proper fractions.
\end{proof}

In toric geometry,  the following fact is well known.
\begin{fact}\upshape{\cite{Fulton}}\label{toriceuler}
Let \(X_\Sigma\) be an \(n\)-dimensional toric variety.
Then $ \chi(X_\Sigma)=\#\Sigma(n)$,
where \(\Sigma(n)\) denotes the set of \(n\)-dimensional cones
in \(\Sigma\).

\end{fact}

By Proposition \ref{main1prop} and Fact \ref{toriceuler}, the following theorem holds.

\begin{thm}\label{formckay}
Let $G$ be a cyclic group of type $\frac{1}{r}(1,a_2,\dots,a_n)$, with $\gcd(r,a_2,\dots,a_n)=1$.
For the Fujiki--Oka resolution $f: Y \to \compl{n}/G$, the following holds
$$
\chi(Y)=h\left(\mathcal{R}_*\left(\frac{(1,a_2,\dots,a_n)}{r}\right)\right) + |G|,
$$
 where $\chi(Y)$ is the topological Euler characteristic of $Y$.
\end{thm}

We apply this theorem to a crepant Fujiki--Oka resolution.

\begin{prop}\cite{SS}\label{prop:ss}
Let $G$ be a cyclic group of type $\frac{1}{r}(1,a_2,\dots,a_n)$ with $\gcd(r,a_2,\dots, a_n)=1$, the Fujiki--Oka resolution of $\compl{n}/G$ is crepant if and only if the ages of all the coefficients of the corresponding remainder polynomial $\mathcal{R}_*\left((1,a_2,...,a_n)/r\right)$
are $1$, where $\age(x)=\sum x_i/r$ for $x=(x_1,\dots,x_n)/r$.
\end{prop}

If $\age(v)$ is $1$, then $h(v)$ equals $0$. Therefore, the following holds.

\begin{cor}
Let $G$ be a cyclic group of type $\frac{1}{r}(1,a_2,\dots,a_n)$ with $\gcd(r,a_2,\dots,a_n)=1$.
If the Fujiki--Oka resolution $f: Y \to \compl{n}/G$ is crepant, then $\chi(Y)=  r$  holds.
\end{cor}

This corollary is the Fujiki--Oka resolution version of  Batyrev's theorem (Theorem \ref{Batyrev}). As shown by the author in \cite{Sato}, there are examples where the Fujiki--Oka resolution satisfies  $\chi(Y)=r$ even if it is not crepant.

\begin{ex}\upshape (\cite{Sato})\label{ex:euler-equality}
Let $k$ be a positive integer and let $G=\frac{1}{6k+1}(1,3,6k-5)$ . Then the Fujiki--Oka resolution of $\compl{3}/G$ is a Hilbert basis resolution, and it has
the Euler characteristic equal to the order of $G$.

\end{ex}

\section{A height criterion for moderate Fujiki--Oka resolutions}
\label{sec:fo-hilbert}
In this section, we characterize the moderateness of
three-dimensional Fujiki--Oka resolutions in terms of the heights
of the coefficients of the remainder polynomial.
We then deduce consequences for Hilbert basis resolutions and
Euler characteristics.

\subsection{Coefficient occurrences and the height criterion}
\label{thm:notation}
Throughout this section, $G$ is a cyclic group of type
$\frac{1}{r}(1,a,b)$ with $\gcd(r,a,b)=1$. 
The chosen generator is
represented by $c_1=(1,a,b)/r$. We work in the lattice
$N_G=\mathbb{Z}^3+\mathbb{Z}c_1$ with the cone
$\sigma_0=\operatorname{Cone}(\bm{e}_1,\bm{e}_2,
\bm{e}_3)$.
We denote by $\mathcal{X}$ the set of monomials with nonzero coefficient
in $\mathcal{R}_*(c_1)$. The empty word \(1\) labels the initial cone.
For \(X\in\mathcal{X}\), the coefficient \(c_X=(1,a_X,b_X)/r_X\) has height \(h_X=h(c_X)=1+a_X+b_X-r_X\).
Distinct monomials are counted separately even when their coefficients
or their associated rays coincide.
In the star subdivision at the Oka center, we obtain the $i$-th \emph{child cone} by replacing the $i$-th generator of the ordered local cone with the center, provided that the resulting cone is three-dimensional.
If this child is singular, the monomial \(Xx_i\) labels it and has coefficient \(c_{Xx_i}=R_i(c_X)\).
Smooth children contribute no terms to the remainder polynomial.
The coefficients of the singular children are
\begin{equation*}
\label{foh:remainder}
\begin{aligned}
 c_{Xx_2}
 &=\frac{1}{a_X}
   \bigl(1,\overline{-r_X}^{\,a_X},
            \overline{b_X}^{\,a_X}\bigr)
 &&(a_X>1),\\
 c_{Xx_3}
 &=\frac{1}{b_X}
   \bigl(1,\overline{a_X}^{\,b_X},
            \overline{-r_X}^{\,b_X}\bigr)
 &&(b_X>1).
\end{aligned}
\end{equation*}
Since taking residues modulo the new denominator preserves the gcd
with that denominator, the remainder formulas give $\gcd(r_{Xx_i},a_{Xx_i},b_{Xx_i})=\gcd(r_X,a_X,b_X)$.
Hence we have \(\gcd(r_X,a_X,b_X)=1\) for every \(X\in\mathcal X\).

For each $X\in\mathcal X$, write the corresponding local cone as
$\operatorname{Cone}(\bm{e}_1,\bm p_X,\bm{q}_X)$ and define its Oka
center by $\displaystyle
 \bm u_X
 =\frac{\bm e_1+a_X\bm p_X+b_X\bm{q}_X}{r_X}$.
The ordered generators are determined recursively by
\[
 (\bm p_1,\bm{q}_1)
   =(\bm e_2,\bm e_3),\ 
 (\bm p_{Xx_2},\bm{q}_{Xx_2})
   =(\bm u_X,\bm{q}_X),\ 
 (\bm p_{Xx_3},\bm{q}_{Xx_3})
   =(\bm p_X,\bm u_X),
\]
where each child formula is used only when $Xx_i\in\mathcal X$.
The triple $(\bm e_1,\bm p_X,\bm{q}_X)$ is a real basis,
but need not be a lattice basis of $N_G$.
For every $X\in\mathcal X$, we have
\[
 N_G=\mathbb Z\bm e_1+\mathbb Z\bm p_X+
      \mathbb Z\bm{q}_X+\mathbb Z\bm u_X,\  \text{and}\ 
 \bigl[N_G:
   \mathbb Z\bm e_1+\mathbb Z\bm p_X+\mathbb Z\bm{q}_X
 \bigr]=r_X.
\]
The three-dimensional cones obtained by subdivision at $\bm u_X$
are
\[
\begin{aligned}
 \sigma_{X,1}
   &=\operatorname{Cone}(\bm u_X,\bm p_X,\bm{q}_X),\\
 \sigma_{X,2}
   &=\operatorname{Cone}(\bm e_1,\bm u_X,\bm{q}_X)
     &&(a_X>0),\\
 \sigma_{X,3}
   &=\operatorname{Cone}(\bm e_1,\bm p_X,\bm u_X)
     &&(b_X>0).
\end{aligned}
\]
The defining relation for the Oka center gives
\(N_G=\mathbb Z\bm{u}_X+\mathbb Z\bm{p}_X+\mathbb Z\bm{q}_X\).
Hence \(\sigma_{X,1}\) is smooth and \(\bm{u}_X\) is primitive.
We call the first child cone \(\sigma_{X,1}\) the \emph{cap} at \(X\).
The other cones, when present, have indices \(a_X\) and \(b_X\).
The denominator strictly decreases at every nonzero remainder
step, so the construction terminates.

Moreover, \(\sigma_{X,1}\) remains a cone of the final fan.
Indeed, every subsequent Oka center in either remaining child has
strictly positive \(\bm{e}_1\)-coordinate in that child's initial
ordered basis. Therefore, it lies outside the face shared with
\(\sigma_{X,1}\), on which this coordinate vanishes.
This argument also applies when one weight is zero and does not
require any height assumption.

Let \(m_{\sigma_{X,1}}\in\operatorname{Hom}(N_G,\mathbb Z)\) be the
linear form taking value one on \(\bm{u}_X,\bm{p}_X,\bm{q}_X\).
Evaluating the Oka center relation gives
\begin{equation*}
  m_{\sigma_{X,1}}(\bm{e}_1)
  =r_X-a_X-b_X
  =1-h_X.
  \label{thm:cap-height}
\end{equation*}

For each smooth three-dimensional cone \(\kappa\), let
\(m_\kappa\in\operatorname{Hom}(N_G,\mathbb Z)\)
take value one on its primitive ray generators.
Their restrictions to the maximal cones of
\(\Sigma_{\mathrm{F.O.R.}}\) define a piecewise linear function
\(\psi:\sigma_0\to\mathbb R\), called the unit function.

Our convention for (upper) convexity is
\[
 \psi(\lambda \bm{v}+(1-\lambda)\bm{v}')
 \ge \lambda\psi(\bm{v})+(1-\lambda)\psi(\bm{v}')
 \qquad
 (\bm{v},\bm{v}'\in\sigma_0,\ 0\le\lambda\le1).
\]
Write \(\Hilb_{N_G}(\sigma_0)\) for the nonzero elements of
\(N_G\cap\sigma_0\) that cannot be expressed as the sum of two nonzero
elements of \(N_G\cap\sigma_0\).

\begin{defi}
\label{thm:moderate}
\upshape
Let \(\Sigma\) be a smooth subdivision of \(\sigma_0\).
Following \cite[Definition~6.6]{CMDY}, we call \(\Sigma\)
\emph{a moderate resolution} of \(\sigma_0\) if the affine
hyperplane $ \{\bm{v}\in N_{G,\mathbb R}\mid m_\kappa(\bm{v})=1\}$
meets every one-dimensional face of \(\sigma_0\), for every
three-dimensional cone \(\kappa\in\Sigma\).
We also call the induced morphism $f:X(N_G,\Sigma)\longrightarrow X(N_G,\sigma_0)$
a moderate toric resolution.
\end{defi}
In the present setting, this condition is equivalent to
\[
 m_\kappa(\bm e_j)>0
 \qquad
 (j=1,2,3)
\]
for every three-dimensional cone \(\kappa\in\Sigma\).
Since these values are integers, it is also equivalent to
\(m_\kappa(\bm e_j)\ge1\) for all such \(\kappa\) and \(j\).

\begin{thm}
\label{thm:criterion}
\upshape
With the notation above, the following conditions are equivalent:
\begin{enumerate}
 \item[(i)]
 \(h_X\le0\) for every \(X\in\mathcal X\), including the empty
 monomial.

 \item[(ii)]
 The unit function \(\psi\) of
 \(\Sigma_{\mathrm{F.O.R.}}\) is convex on \(\sigma_0\).

 \item[(iii)]
 The Fujiki--Oka resolution associated with \(c_1\) is moderate.
\end{enumerate}
\end{thm}

\begin{proof}[The implications
  \(\mathrm{(ii)}\Rightarrow\mathrm{(iii)}
  \Rightarrow\mathrm{(i)}\)]
The original generators \(\bm{e}_1,\bm{e}_2,\bm{e}_3\) are primitive and remain
rays of the final fan, so \(\psi(\bm{e}_j)=1\).
If \(\psi\) is convex, then we have $\psi=\min_{\kappa}m_{\kappa},$
where \(\kappa\) ranges over the maximal cones.
Consequently,
\[
  m_{\kappa}(\bm{e}_j)\ge\psi(\bm{e}_j)=1
  \qquad (j=1,2,3),
\]
which proves moderateness.
Next, suppose that the resolution is moderate. Then every cap
\(\sigma_{X,1}\) satisfies
\[
  0<m_{\sigma_{X,1}}(\bm{e}_1)=1-h_X.
\]
Since \(h_X\) is an integer, we obtain \(h_X\le0\) for every
\(X\in\mathcal X\).
\end{proof}

To prove the remaining implication, we first compare the height of
a coefficient with that of a singular child.

\begin{lem}
\label{lem:height-lemma}
\upshape
Suppose that \(X,Xx_2\in\mathcal X\) and
\(h_X,h_{Xx_2}\le0\). Then $-h_{Xx_2}=\overline{-h_X}^{\,a_{X}}$ and $h_X\le h_{Xx_2}\le0$.
The analogous statement holds for $x_3$, with $b_X$
in place of $a_X$.
\end{lem}

\begin{prf}
By the definition of remainder map, 
\[
 h_{Xx_2}
 =1+\overline{-r_X}^{\,a_X}+\overline b_X^{\,a_X}-a_X
 \equiv h_X\pmod {a_X}.
\]
The residues are nonnegative, and \(h_{Xx_2}\le0\) by assumption.
Hence we have $ 1-a_X\le h_{Xx_2}\le0$.
There is a unique integer in this interval congruent to \(h_X\)
modulo \(a_X\), so $ h_{Xx_2}=-\overline{-h_X}^{\,a_X}.$
Since \(h_X\le0\), we have $ 0\le\overline{-h_X}^{\,a_X}\le-h_X$,
which gives \(h_X\le h_{Xx_2}\le0\).
The proof for \(x_3\) is identical.
\end{prf}

It remains to prove (i) \(\Rightarrow\) (ii) in Theorem~\ref{thm:criterion}.
We argue by induction on the local index.
Convexity within each child cone follows from the induction
hypothesis, with smooth children requiring no further argument.
We need only compare the unit functions across the
walls separating the children.
Lemma~\ref{lem:height-lemma} controls the walls adjacent to the cap, while
Lemma~\ref{lem:common-face-convexity} treats the common face of the other two children.

\subsection{Convexity across the common face}

To study convexity across the common face of the two child cones,
we compare the two boundary branches with those of an auxiliary
type whose second and third weights are coprime.

Fix \(X\in\mathcal{X}\) with \(a_X,b_X>0\) and \(h_X\le0\). In a calculation at a fixed monomial \(X\), for simplicity, we omit the subscript $X$. 
The child cones \(\sigma_{X,2}\) and \(\sigma_{X,3}\) meet along
\(F=\operatorname{Cone}(\bm{e}_1,\bm{u})\).
To compare their subdivisions along \(F\), we write \(a=g\alpha\) and
\(b=g\beta\), where \(g=\gcd(a,b)\). Then \(\gcd(\alpha,\beta)=1\).
Write \(r=g\rho-\delta\) with \(0\le\delta<g\), so that
\(\rho=\lceil r/g\rceil\), and set \(\bm{t}=(\bm{e}_1+\delta \bm{u})/g\).
Substituting \(\bm{e}_1=g\bm{t}-\delta \bm{u}\) into \(r\bm{u}=\bm{e}_1+a\bm{p}+b\bm{q}\) gives
\[
\rho \bm{u}=\bm{t}+\alpha \bm{p}+\beta \bm{q}.
\]
Thus \(\bm{t}=\rho \bm{u}-\alpha \bm{p}-\beta \bm{q}\) lies in \(N_G\).
Since \(h\le0\), we also have \(\rho\ge1+\alpha+\beta\).
In the lattice basis \((\bm{u},\bm{p},\bm{q})\), the last two coordinates of a
lattice point in \(\operatorname{span}_{\mathbb R}F\) are proportional
to \((\alpha,\beta)\).
Since \(\gcd(\alpha,\beta)=1\), they are an integral multiple
of this pair. Hence
\begin{equation*}
 N_F:=N_G\cap\operatorname{span}_{\mathbb R}F
 =\mathbb Z\bm{t}+\mathbb Z\bm{u}.
 \label{eq:boundary-face-lattice}
\end{equation*}
Thus \(\bm{t}\) is primitive. In the lattice \(N_G\), the three-dimensional cone \(\operatorname{Cone}(\bm{t},\bm p,\bm{q})\) defines a cyclic quotient
singularity of type \(\frac{1}{\rho}(1,\alpha,\beta)\).
We call this the \emph{auxiliary type}.
With \(\bm{t}\) as the distinguished generator, the Oka center of
this cone is again \(\bm u\), since
\(\rho\bm u=\bm{t}+\alpha\bm p+\beta\bm{q}\).
When \(g>1\), the ray \(\mathbb R_{\geq0}\bm{t}\) is adjacent to
\(\mathbb R_{\geq0}\bm{u}\) in the minimal subdivision of \(F\),
as shown in Figure~\ref{fig:common-face}.

\begin{rem}
\label{thm:surface-facts}
\upshape
Let $ F=\operatorname{Cone}(\bm{w}_1,\bm{w}_2)$ and $N_F=\mathbb Z\bm{w}_1+\mathbb Z\bm{w}_2      +\mathbb Z\frac{\bm{w}_1+s\bm{w}_2}{g}$,
where \(0<s<g\) and \(\gcd(g,s)=1\).
By the Hirzebruch--Jung construction and the convex-hull
description of the minimal resolution
\cite[Theorems~10.2.3 and~10.2.8]{CLS},
its unit function \(\psi_F\) is convex.
Every linear piece \(\ell\) is integral on \(N_F\) and satisfies
\[
 1\le\ell(\bm{w}_1)\le g-s,
 \qquad
 \ell(\bm{w}_2)\ge1.
 \label{thm:surface-bounds}
\]
The unit form on the cone adjacent to \(\bm{w}_2\), generated by
\(\bm{w}_2\) and \((\bm{w}_1+s\bm{w}_2)/g\), takes value \(g-s\) at \(\bm{w}_1\).
By convexity, the slopes along \(\bm{w}_2+t\bm{w}_1\), \(t\ge0\),
are nonincreasing, which gives \(\ell(\bm{w}_1)\le g-s\).
\end{rem}

While the Oka center lies outside \(F\), only the third child
on the \(\sigma_{X,2}\) side and the second child on the
\(\sigma_{X,3}\) side retain the whole face \(F\).
Thus the branches relevant to the wall comparison are
\(Xx_2x_3^j\) and \(Xx_3x_2^j\).

\begin{lem}
\label{lem:boundary-reduction}
\upshape
Suppose that every nonzero coefficient of $\mathcal{R}_*(c_X)$,
including its constant term $c_X$, has nonpositive height.
Then every nonzero coefficient on the boundary branches $x_2x_3^j$
and $x_3x_2^j$ $(j\geq 0)$ of the auxiliary type $(1,\alpha,\beta)/\rho$
has nonpositive height.

Each nonzero term on these auxiliary branches has the same Oka center
as the original term whose monomial is obtained by adding the prefix
\(X\). On each side of \(F\), the original and auxiliary branches
also have the same terminal outside generator.
\end{lem}

\begin{prf}
On the auxiliary branch \(x_2x_3^j\), a nonzero coefficient
\((1,A,B)/Q\) corresponds to the original coefficient
\[
 \frac{(1,gA+\delta,gB)}{gQ}.
\]
This holds at \(x_2\) because \(r=g\rho-\delta\), and is preserved
at every nonzero auxiliary \(R_3\)-step, since
\[
 \overline{gA+\delta}^{\,gB}
 =g\overline A^{\,B}+\delta,
 \qquad
 \overline{-gQ}^{\,gB}=g\overline{-Q}^{\,B}.
\]
If \(\bm{T}\) is the current outside generator, the lifted centers agree:
\[
 \frac{\bm{e}_1+(gA+\delta)\bm{u}+gB\bm{T}}{gQ}
 =\frac{\bm{t}+A\bm{u}+B\bm{T}}{Q}.
\]
Moreover, \(1\leq1+\delta\leq g\), so
\[
 1+gA+\delta+gB\leq gQ
 \quad\Longleftrightarrow\quad
 1+A+B\leq Q.
\]
Thus the corresponding auxiliary coefficient has nonpositive height.
On the other branch the original coefficient is
\((1,gA,gB+\delta)/(gQ)\), and the argument is symmetric.

Since \(\gcd(\alpha,\beta)=1\), the auxiliary branches end in
smooth cones
\[
 \operatorname{Cone}(\bm{t},\bm{u},\bm{V}_2),\qquad
 \operatorname{Cone}(\bm{t},\bm{u},\bm{V}_3).
\]
The denominators and active weights on the original branches are
those on the auxiliary branches multiplied by \(g\).
When \(g>1\), the original branches therefore reach coefficients
\((1,\delta,0)/g\) and \((1,0,\delta)/g\), respectively, with the
same outside generators \(\bm{V}_2,\bm{V}_3\).
All subsequent Oka centers in these terminal cones lie in \(F\).
When \(g=1\), the original branches already end in the same smooth cones.
If \(\alpha=1\) or \(\beta=1\), the corresponding auxiliary child
is smooth from the outset, and the same conclusion holds.
\end{prf}

\begin{rem}
 \label{thm:wall-test}
\upshape
For a smooth fan with full-dimensional convex support,
convexity of the unit function can be checked at the
interior walls; see \cite[Lemma~6.1.5]{CLS}.

For adjacent smooth cones
\(\operatorname{Cone}(\bm{V}_2,\bm{P}_1,\bm{P}_2)\) and
\(\operatorname{Cone}(\bm{V}_3,\bm{P}_1,\bm{P}_2)\), write
\(\bm{V}_2+\bm{V}_3=\lambda \bm{P}_1+\mu \bm{P}_2\), where \(\lambda,\mu\in\mathbb Z\).
Let \(\ell\) be the linear form on their common wall satisfying
\(\ell(\bm{P}_1)=\ell(\bm{P}_2)=1\).
The unit function is convex across this wall if and only if
\[
 \ell(\bm{V}_2+\bm{V}_3)=\lambda+\mu\ge2.
\]
Indeed, the unit form \(m\) of the first cone satisfies
\(m(\bm V_3)=\lambda+\mu-1\).
\end{rem}

\begin{lem}
\label{lem:common-face-convexity}
\upshape
Under the assumptions of Lemma~\ref{lem:boundary-reduction},
the final unit functions on \(\sigma_{X,2}\) and \(\sigma_{X,3}\)
glue convexly across every wall contained in
\(F\).
\end{lem}

\begin{prf}
Let \(\bm{V}_2\) and \(\bm{V}_3\) be the outside generators given by
Lemma~\ref{lem:boundary-reduction}. The vectors \(\bm{t}\) and
\(\bm{u}\) form a basis of \(N_F\). For \(i=2,3\), the cone
\(\operatorname{Cone}(\bm{t},\bm{u},\bm{V}_i)\) in the auxiliary
subdivision is smooth, so \(\bm{V}_i\) extends this basis to a basis
of \(N_G\).

For \(g>1\), the local types at the ends of the two original branches
are \((1,\delta,0)/g\) and \((1,0,\delta)/g\). In each case the zero
weight corresponds to the outside generator. Thus both branches induce
on \(F\) the fan of the minimal resolution of the surface type
\((1,\delta)/g\), where \(\gcd(g,\delta)=\gcd(g,r)=1\).
For \(g=1\), the face \(F\) is smooth and remains unsubdivided.

In either case, each maximal cone \(\tau\) of the surface fan is the
common wall of \(\operatorname{Cone}(\tau,\bm{V}_2)\) and
\(\operatorname{Cone}(\tau,\bm{V}_3)\).
The unit form \(\ell_\tau\) on \(\operatorname{span}_{\mathbb R}F\)
takes the value \(1\) on both primitive ray generators of \(\tau\).
By Remark~\ref{thm:wall-test}, it suffices to show that
\(\ell_\tau(\bm{V}_2+\bm{V}_3)\geq 2\) for every \(\tau\).

We begin with \(\tau=\operatorname{Cone}(\bm{t},\bm{u})\).
For \(g>1\), this is the lower wall in Figure~\ref{fig:common-face}.
For \(g=1\), it is \(F\) itself.
The auxiliary cones \(\operatorname{Cone}(\bm{t},\bm{u},\bm{q})\)
and \(\operatorname{Cone}(\bm{t},\bm{p},\bm{u})\) have indices
\(\alpha\) and \(\beta\), respectively.
Since $ [N_G:\mathbb Z\bm{t}+\mathbb Z\bm{u}+\mathbb Z\bm{q}]=\alpha$,
we have
\(\alpha\bm{V}_2\in
\mathbb Z\bm{t}+\mathbb Z\bm{u}+\mathbb Z\bm{q}\).
Moreover, \(\bm{V}_2\in
\operatorname{Cone}(\bm{t},\bm{u},\bm{q})\), so we may write $
\alpha\bm{V}_2
 =k_2\bm{t}+n_2\bm{u}+d_2\bm{q}$,
for nonnegative integers $k_2,n_2,d_2$.
Computing determinants with respect to a lattice basis of \(N_G\), the smoothness of
\(\operatorname{Cone}(\bm{t},\bm{u},\bm{V}_2)\) gives $\displaystyle
 1
 =\left|\det(\bm{t},\bm{u},\bm{V}_2)\right|
 =\frac{d_2}{\alpha}
   \left|\det(\bm{t},\bm{u},\bm{q})\right|
 =d_2$.
Thus, applying the same argument on the other side, we obtain
\[
 \alpha\bm{V}_2=k_2\bm{t}+n_2\bm{u}+\bm{q},
 \qquad
 \beta\bm{V}_3=k_3\bm{t}+n_3\bm{u}+\bm{p},
 \qquad k_i,n_i\in\mathbb Z_{\ge0}\quad(i=2,3).
\]
Using \(\alpha \bm{p}+\beta \bm{q}=\rho \bm{u}-\bm{t}\), we have
\[
  \bm{V}_2+\bm{V}_3=\lambda \bm{t}+\mu \bm{u},\qquad
  \lambda=\frac{\beta k_2+\alpha k_3-1}{\alpha\beta},
  \qquad
  \mu=\frac{\rho+\beta n_2+\alpha n_3}{\alpha\beta}.
\]
Since \((\bm{t},\bm{u})\) is a basis of \(N_F\), both \(\lambda\) and \(\mu\)
are integers. In particular, \(\mu\ge1\).

We claim that \(\lambda+\mu\ge2\).
If \(\alpha,\beta>1\), then \(\beta k_2+\alpha k_3\) is nonnegative
and congruent to \(1\) modulo \(\alpha\beta\). Since
\(\alpha,\beta>1\), this sum cannot equal \(1\), so it is at least
\(\alpha\beta+1\). Hence \(\lambda\ge1\).
If \(\min(\alpha,\beta)=1\), then \(h\le0\) gives
\(\rho\ge1+\alpha+\beta\), and hence we have
\[
  \lambda+\mu
  \ge \frac{\rho-1}{\alpha\beta}
  \ge \frac{\alpha+\beta}{\alpha\beta}
  = 1+\frac{1}{\alpha\beta}>1.
\]
The claim follows by integrality.

Now let \(\tau\) be any maximal cone of the surface fan.
For \(g>1\), Remark~\ref{thm:surface-facts} gives
\(\ell_\tau(\bm{e}_1)\le g-\delta\) and \(\ell_\tau(\bm{u})\ge1\).
Since \(\ell_\tau\) takes positive integral values on the nonzero
lattice points of \(F\), we obtain
\[
  1\le\ell_\tau(\bm{t})
  =\frac{\ell_\tau(\bm{e}_1)+\delta\ell_\tau(\bm{u})}{g}
  \le\frac{g-\delta+\delta\ell_\tau(\bm{u})}{g}
  \le\ell_\tau(\bm{u}).
\]
For \(g=1\), both values are one. Therefore, we have
\[
  \ell_\tau(\bm{V}_2+\bm{V}_3)
  =(\lambda+\mu)\ell_\tau(\bm{t})
   +\mu\bigl(\ell_\tau(\bm{u})-\ell_\tau(\bm{t})\bigr)
  \ge2.
\]
This proves convexity across every wall in \(F\).
\end{prf}

\begin{figure}
  \centering
\begin{tikzpicture}[
  x=1cm,y=1cm,font=\small,
  line cap=round,line join=round,
  fo mesh/.style={draw=foMesh,line width=.35pt},
  fo outline/.style={draw=black!75,line width=.65pt},
  fo wall/.style={draw=foRed,line width=1.3pt},
  fo label/.style={fill=white,inner sep=1.5pt},
  fo dot/.style={circle,inner sep=0pt,minimum size=3.1pt}
]

  \begin{scope}[xshift=6.2cm]
    \coordinate (E1) at (0,4);
    \coordinate (E2) at (-3,0);
    \coordinate (E3) at (3,0);
    \coordinate (U) at ({3*(10-6)/17},{4*1/17}); 
    \coordinate (A) at ({3*(13-1)/17},{4*3/17}); 
    \coordinate (V2) at ({3*(9-2)/17},{4*6/17}); 
    \coordinate (Fpoint) at ({3*(5-3)/17},{4*9/17}); 
    \coordinate (B) at ({3*(3-12)/17},{4*2/17}); 
    \coordinate (V3) at ({3*(6-7)/17},{4*4/17}); 
    \coordinate (C) at ({3*(2-8)/17},{4*7/17}); 
    \coordinate (D) at ({3*(1-4)/17},{4*12/17}); 
    \fill[foBlue!7] (E1)--(Fpoint)--(V2)--cycle;
    \fill[foOrange!7] (E1)--(Fpoint)--(V3)--cycle;
    \fill[foBlue!22] (Fpoint)--(U)--(V2)--cycle;
    \fill[foOrange!24] (Fpoint)--(U)--(V3)--cycle;
    \draw[fo mesh] (A) -- (E1);
    \draw[fo mesh] (A) -- (E3);
    \draw[fo mesh] (A) -- (U);
    \draw[fo mesh] (A) -- (V2);
    \draw[fo mesh] (B) -- (C);
    \draw[fo mesh] (B) -- (E2);
    \draw[fo mesh] (B) -- (U);
    \draw[fo mesh] (B) -- (V3);
    \draw[fo mesh] (C) -- (D);
    \draw[fo mesh] (C) -- (E2);
    \draw[fo mesh] (C) -- (V3);
    \draw[fo mesh] (D) -- (E1);
    \draw[fo mesh] (D) -- (E2);
    \draw[fo mesh] (D) -- (V3);
    \draw[fo mesh] (E1) -- (E2);
    \draw[fo mesh] (E1) -- (E3);
    \draw[fo mesh] (E1) -- (Fpoint);
    \draw[fo mesh] (E1) -- (V2);
    \draw[fo mesh] (E1) -- (V3);
    \draw[fo mesh] (E2) -- (E3);
    \draw[fo mesh] (E2) -- (U);
    \draw[fo mesh] (E3) -- (U);
    \draw[fo mesh] (Fpoint) -- (U);
    \draw[fo mesh] (Fpoint) -- (V2);
    \draw[fo mesh] (Fpoint) -- (V3);
    \draw[fo mesh] (U) -- (V2);
    \draw[fo mesh] (U) -- (V3);
    \draw[fo outline] (E1)--(E2)--(E3)--cycle;
    \draw[foBlue,line width=.7pt] (Fpoint)--(V2)--(U);
    \draw[foOrange,line width=.7pt] (Fpoint)--(V3)--(U);
    \draw[foRed,line width=.65pt] (E1)--(Fpoint);
    \draw[fo wall] (Fpoint)--(U);
    \node[text=foBlue,inner sep=0pt]
      at (barycentric cs:Fpoint=1,U=1,V2=1) {};
    \node[text=foOrange,inner sep=0pt]
      at (barycentric cs:Fpoint=1,U=1,V3=1) {};
    \foreach \p in {E1,E2,E3,A,B,C,D}
      \fill[black!55] (\p) circle[radius=1.15pt];
    \foreach \p in {Fpoint,U}
      \node[fo dot,fill=foRed,draw=white,line width=.3pt] at (\p) {};
    \node[fo dot,fill=foBlue,draw=white,line width=.3pt] at (V2) {};
    \node[fo dot,fill=foOrange,draw=white,line width=.3pt] at (V3) {};
    \node[above=4pt] at (E1) {$\bm{e}_1$};
    \node[below=5pt] at (E2) {$\bm{p}$};
    \node[below=5pt] at (E3) {$\bm{q}$};
    \draw[foRed,line width=.45pt] (Fpoint)--(-.57,2.40);
    \node[fo label,text=foRed,anchor=east] at (-.57,2.40) {$\bm{t}$};
    \draw[foRed,line width=.45pt] (U)--(.86,-.20);
    \node[below=1pt,text=foRed] at (.86,-.20) {$\bm{u}$};
    \draw[foBlue,line width=.45pt] (V2)--(1.74,1.61);
    \node[fo label,text=foBlue,anchor=west] at (1.74,1.61) {$\bm{V}_2$};
    \draw[foOrange,line width=.45pt] (V3)--(-.85,.81);
    \node[fo label,text=foOrange,anchor=east] at (-.85,.81) {$\bm{V}_3$};
    \draw[foRed,line width=.45pt] ($(E1)!.25!(U)$)--(1.10,3.20);
    \node[fo label,text=foRed,anchor=west] at (1.10,3.20) {$F$};
  \end{scope}
 

\end{tikzpicture}
  \caption{Adjacent cones sharing the wall
\(\operatorname{Cone}(\bm{t},\bm{u})\subset F\).}
  \label{fig:common-face}
\end{figure}

\subsection{Proof of the coefficient-height criterion}
\label{thm:criterion-proof}
\begin{proof}[Proof of
  \(\mathrm{(i)}\Rightarrow\mathrm{(ii)}\)
  in Theorem~\ref{thm:criterion}]
Assume \(\mathrm{(i)}\).
We prove by induction on \(r_X\) that the final unit function
restricted to
\(\operatorname{Cone}(\bm{e}_1,\bm{p},\bm{q})\) is convex whenever
\(c_X\) and all its nonzero descendants have nonpositive height.
We use the local notation \(\left(r,a,b,h,\bm{p},\bm{q},\bm{u} \right)\) at \(X\).

If \(a=0\) or \(b=0\), the remaining surface type is coprime,
since \(\gcd(r,a,b)=1\).
The Fujiki--Oka construction is the minimal surface resolution
times the unused axis.
Its unit function is the sum of the surface unit function and
the coordinate along that axis, and is therefore convex.
Suppose that \(a,b>0\).
The first subdivision consists of
\(\sigma_{X,1},\sigma_{X,2},\sigma_{X,3}\).

The cap \(\sigma_{X,1}\) is a smooth cone of the final fan.
The other two cones have indices \(a,b<r\), so their final unit
functions are convex by induction; a child of index one is
already smooth.
It remains to check the walls between these three regions.

Let \(m'\) be the unit form of the final cone across
\(\operatorname{Cone}(\bm{u},\bm{q})\) from the cap. Then \(m'(\bm{u})=m'(\bm{q})=1\).
If \(a>1\), this cone is the cap at \(Xx_2\), so
\(m'(\bm{e}_1)=1-h_{Xx_2}\). By Lemma~4.3, \(h\le h_{Xx_2}\le0\), and hence we have $1\le m'(\bm{e}_1)\le1-h$.

If \(a=1\), the neighboring cone is \(\sigma_{X,2}\).
It contains \(\bm{e}_1\) as a primitive ray generator, so \(m'(\bm{e}_1)=1\)
and the same bounds hold.
Evaluating \(r\bm{u}=\bm{e}_1+a\bm{p}+b\bm{q}\), we obtain
\[
  m'(\bm{p})
  =\frac{r-b-m'(\bm{e}_1)}a
  =1+\frac{1-h-m'(\bm{e}_1)}a
  \ge1.
\]
The wall criterion proves convexity across
\(\operatorname{Cone}(\bm{u},\bm{q})\).
The argument across \(\operatorname{Cone}(\bm{p},\bm{u})\) is symmetric.
Every remaining wall between \(\sigma_{X,2}\) and \(\sigma_{X,3}\)
is contained in \(F=\operatorname{Cone}(\bm{e}_1,\bm{u})\).
Every nonzero coefficient of $\mathcal{R}_*(c_X)$ is of the form
$c_{XW}$ for some word $W$ or $1$. By~(i), each such
coefficient has nonpositive height. Therefore, Lemma~\ref{lem:common-face-convexity} gives
convexity across every wall contained in $F$.
Thus the unit function is convex across every interior wall,
which completes the induction.
Taking \(X=1\) proves convexity on \(\sigma_0\).
\end{proof}

\begin{cor}
\label{thm:hilbert-corollary}
\upshape
If the equivalent conditions of Theorem~\ref{thm:criterion}
hold, then the primitive ray generators of
\(\Sigma_{\mathrm{F.O.R.}}\) are exactly
\(\Hilb_{N_G}(\sigma_0)\).
In particular, the Fujiki--Oka resolution is a Hilbert basis
resolution.
\end{cor}

\begin{prf}
By Theorem~\ref{thm:criterion}, the resolution is moderate.
The conclusion follows from
\cite[Proposition~6.7]{CMDY}, which states that every moderate
toric resolution is a Hilbert basis resolution.
\end{prf}
\begin{rem}
\label{thm:euler-corollary}
\upshape
Let \(f:Y\to\mathbb C^3/G\) be the Fujiki--Oka resolution.
If \(f\) is moderate, Theorems~\ref{formckay}
and~\ref{thm:criterion} give $ \chi(Y)=r+\sum_{X\in\mathcal X}h_X\le r$.
Among moderate Fujiki--Oka resolutions, by Proposition \ref{prop:ss}, equality holds
if and only if the resolution is crepant.
Without the moderateness assumption, however,
\(\chi(Y)=r\) does not imply that the resolution is crepant (see Example~\ref{ex:euler-equality}).
\end{rem}

\begin{ex}\upshape
  Let \(G=\frac{1}{17}(1,6,9)\). At the initial cone we have $\bm{u}=\frac{1}{17}(1,6,9)$ and $(g,\alpha,\beta,\rho,\delta)=(3,2,3,6,1)$.
The monomials \(x_2x_3\) and \(x_3x_2^2\) give the same Oka center $\bm{t}=\dfrac{\bm{e}_1+\bm{u}}{3}=\dfrac{1}{17}(6,2,3)$.
Since \((\bm{t},\bm{u})\) is a basis of \(N_F\) and \(\bm{e}_1=3\bm{t}-\bm{u}\), the
rays of the minimal surface fan are generated by \(\bm{u},\bm{t},\bm{e}_1\). The two boundary branches give the outside generators $\bm{V}_2=\dfrac{\bm{t}+\bm{e}_3}{2}=\dfrac{1}{17}(3,1,10)$ and $\bm{V}_3=\dfrac{2\bm{t}+\bm{e}_2}{3}=\dfrac{1}{17}(4,7,2)$.
In particular, we have \(\bm{V}_2+\bm{V}_3=\bm{t}+\bm{u}\).

Let \(\ell_1\) and \(\ell_2\) be the unit forms of
\(\operatorname{Cone}(\bm{u},\bm{t})\) and \(\operatorname{Cone}(\bm{t},\bm{e}_1)\),
respectively. We have \(\ell_1(\bm{u})=\ell_1(\bm{t})=1\), while
\(\ell_2(\bm{t})=\ell_2(\bm{e}_1)=1\) gives
\(\ell_2(\bm{u})=\ell_2(3\bm{t}-\bm{e}_1)=2\). Consequently, 
\[
  \ell_1(\bm{V}_2+\bm{V}_3)=2,\qquad \ell_2(\bm{V}_2+\bm{V}_3)=3.
\]
Thus the unit function is linear across the first wall and bends
strictly across the second.
The six coefficients of the remainder polynomial all have height \(-1\).
By Theorem~\ref{thm:criterion}, the Fujiki--Oka resolution \(Y\) is moderate.
The Euler characteristic formula then gives
\(\chi(Y)=17+6(-1)=11\).
\end{ex}

When \(g>1\), Remark \ref{thm:surface-facts} identifies \(\bm{t}\) as the primitive generator
of the ray adjacent to \(\bm{u}\) in the minimal surface fan of \(F\).

\section{Examples and remarks}

The following examples illustrate the distinctions between the
Hilbert basis property, moderateness, and the existence of a moderate
resolution. For each presentation $\frac1r(1,a,b)$, we take the first
coordinate axis as the distinguished axis.
We recall that a moderate resolution cannot exist if
$\operatorname{Hilb}_{N_G}(\sigma_0)$ contains an element whose coordinate sum exceeds one. 

\paragraph{(i) A terminal example.}
Let $G=\frac12(1,1,1)$ and $\bm{v}=\frac12(1,1,1)$.
The Fujiki--Oka resolution is the star subdivision at $\bm{v}$, and $\operatorname{Hilb}_{N_G}(\sigma_0)=\{\bm{e}_1,\bm{e}_2,\bm{e}_3,\bm{v}\}$.
Thus it is a Hilbert basis resolution, whereas no moderate
resolution exists since $\operatorname{age}(\bm{v})=3/2$.

\paragraph{(ii) A non-canonical example.}
We next consider $G=\frac15(1,1,2)$ and write
$\bm{v}_1=\frac15(1,1,2)$ and $\bm{v}_3=\frac15(3,3,1)$.
We have
\[
 \mathcal R_*\!\left(\frac{(1,1,2)}5\right)
 =\frac{(1,1,2)}5+\frac{(1,1,1)}2x_3.
\]
The corresponding Oka centers are $\bm{v}_1$ and
$\bm{v}_3=(\bm{e}_1+\bm{e}_2+\bm{v}_1)/2$, and the Hilbert basis is
$\{\bm{e}_1,\bm{e}_2,\bm{e}_3, \bm{v}_1,\bm{v}_3\}$.
Hence the Fujiki--Oka resolution is a Hilbert basis resolution.
Since $\operatorname{age}(\bm{v}_1)=4/5$ and
$\operatorname{age}(\bm{v}_3)=7/5$, the singularity is non-canonical
and admits no moderate resolution.
Nevertheless, the coefficient heights are $-1$ and $1$, so
$\chi(Y)=5=|G|$.

\paragraph{(iii) Existence of a different moderate resolution.}
For $G=\frac1{14}(1,4,7)$, we use the lattice points $
 \bm{v}_1=\dfrac{(1,4,7)}{14},
 \bm{v}_2=\dfrac{(2,8,0)}{14},
 \bm{v}_4=\dfrac{(4,2,0)}{14},$ and 
 $\bm{v}_7=\dfrac{(7,0,7)}{14}$.
The Fujiki--Oka center associated with $x_2$ is
$(\bm{e}_1+2\bm{v}_1+3\bm{e}_3)/4=\bm{v}_4+\bm{e}_3$.
Thus this resolution is not a Hilbert basis resolution and is
not moderate.
On the other hand, successive star subdivisions at $\bm{v}_4,\bm{v}_2,\bm{v}_1,\bm{v}_7$
give a moderate resolution, as we can check directly from
its maximal cones:
\[
\begin{gathered}
 \operatorname{Cone}(\bm{e}_1,\bm{v}_7,\bm{v}_4),\quad
 \operatorname{Cone}(\bm{v}_7,\bm{e}_3,\bm{v}_4),\quad
 \operatorname{Cone}(\bm{e}_2,\bm{e}_3,\bm{v}_1),\\
 \operatorname{Cone}(\bm{e}_3,\bm{v}_4,\bm{v}_1),\quad
 \operatorname{Cone}(\bm{v}_4,\bm{v}_2,\bm{v}_1),\quad
 \operatorname{Cone}(\bm{v}_2,\bm{e}_2,\bm{v}_1).
\end{gathered}
\]
\paragraph{(iv) Dependence on the distinguished axis.}
 Let $G=\frac1{19}(1,5,12)$.
The coefficient at $x_2x_3$ is $(1,1,1)/2$, of height one,
so the corresponding Fujiki--Oka resolution is not moderate.
Replacing the group generator by its fourth power and
interchanging the first two coordinates gives the presentation
$\frac1{19}(1,4,10)$.
All nonzero coefficients of its remainder polynomial have
height zero or $-4$. By Theorem~\ref{thm:criterion}, the corresponding
Fujiki--Oka resolution is moderate.

Yusuke Sato\\
\sc{Department of Mathematics, Osaka Institute of Technology, 5-16-1 Ohmiya,
Asahi-ku, Osaka, 535-8585, Japan}\\
Email address: yusuke.sato@oit.ac.jp

\end{document}